\pdfoutput=1
\documentclass{amsart}

\usepackage[utf8]{inputenc}
\usepackage[T1]{fontenc}
\usepackage{lmodern}
\usepackage{amssymb}
\usepackage{mathtools}
\usepackage{lscape}
\usepackage{cite}
\usepackage{hyperref}

\newtheorem{theorem}{Theorem}[section]

\theoremstyle{definition}
\newtheorem{definition}[theorem]{Definition}

\newtheorem{example}[theorem]{Example}

\newtheorem{question}[theorem]{Question}

\theoremstyle{remark}
\newtheorem{remark}[theorem]{Remark}

\numberwithin{equation}{section}

\begin{document}

\title{Complete intersection monomial curves and free resolutions}

\author{Feza Arslan}
\address{Department of Mathematics,
Mimar Sinan Fine Arts University, Istanbul, 34349 Turkey}
\address{Current Address: Department of Mathematics,
Koç University, Istanbul, 34450 Turkey}
\email{sefaarslan@ku.edu.tr, sefafezaarslan@gmail.com}
\author{Özgür İnce}
\address{Department of Mathematics, Sivas Cumhuriyet University, Sivas, 58050 Turkey}
\email{ozgince@cumhuriyet.edu.tr, ozgince@gmail.com}

\thanks{This is an Accepted Manuscript of an article published by
Taylor \& Francis in \emph{Communications in Algebra} on 26 November 2025,
available at: \url{https://doi.org/10.1080/00927872.2025.2583313}.
Published version: Comm.\ Algebra \textbf{54}(6) (2026), 2432--2441.
This version is deposited under the terms of the Creative Commons
Attribution-NonCommercial-NoDerivatives License
(\url{http://creativecommons.org/licenses/by-nc-nd/4.0/}). This version corrects some minor typographical errors in the published version.}




\commby{}

\dedicatory{}

\begin{abstract} In this article, we compute the minimal free resolutions of the tangent cones of complete intersection monomial curves in affine 4-space under certain conditions. We first show that a complete intersection monomial curve in 3-space has a complete intersection tangent cone if and only if it is obtained by a nice gluing. By using the technique of nice gluing, we construct infinitely many families of complete intersection monomial curves with complete intersection tangent cones. We also construct infinitely many families of the complete intersection monomial curves that do not have complete intersection tangent cones and compute the minimal free resolutions of these families.  
\end{abstract}

\maketitle

\section{Introduction}

In this article, our aim is to compute the minimal free resolutions of the tangent cones of complete intersection monomial curves in affine 4-space. A monomial curve $C$ in affine $l$-space over $K$ is a curve having parametrization 
 $$x_1=t^{n_1}, x_2=t^{n_2},\ldots, x_l=t^{n_l}$$
 for some positive integers $n_1 < \cdots < n_l$ with $gcd(n_1,\ldots,n_l)=1$. The corresponding local ring is $K[[t^{n_1}, t^{n_2},\ldots, t^{n_l}]]$ and the corresponding numerical semigroup $S$ is
$$S=\displaystyle \langle n_{1},n_{2},\ldots,n_{l}\rangle=
                                                        \left \{ \sum_{i=1}^l a_{i}n_{i}, \; a_{i}\in\mathbb{Z}_{\geq0} \right \}$$
                                                         
Let $I(C) \subset K[x_{1}, x_2, \ldots, x_l]$ denote the defining ideal of the monomial curve $C$ and $f_*$ denote the homogeneous component of $f \in K[x_{1}, x_2, \ldots, x_l]$ of the lowest degree.
Then the homogeneous ideal $I(C)_*$ formed by $f_*$ for all $f\in I(C)$ is called the defining ideal of the tangent cone of $C$. It is well-known that the coordinate ring of the tangent cone $K[x_1, x_2,\ldots , x_l]/I(C)_*$ is isomorphic to the associated graded ring of the corresponding local ring $K[[t^{n_1}, t^{n_2}, \ldots, t^{n_l}]]$. Recall that the associated graded ring of a local ring $R$ with maximal ideal $\mathfrak{m}$, denoted by $gr_{\mathfrak{m}}(R)$ is given as
$\displaystyle gr_{\mathfrak{m}}(R) = \bigoplus_{i=0}^{\infty} \mathfrak{m}^i/\mathfrak{m}^{i+1}$.
 
For a monomial curve $C$, there is a rich literature on the minimal free resolutions of both $K[x_1, x_2,\ldots , x_l]/I(C)$ and $K[x_1, x_2,\ldots , x_l]/I(C)_*$. Stamate gives an excellent survey of this rich literature \cite{stamate}. Various authors contributed to this literature \cite{barucci}, \cite{enescu, gimenez1, gimenez},
\cite{herzog1,jafari,jayanthan,kumar,paper21,paper22,oneto}, \cite{saha, sharifan, sharifan1}, \cite{ss, sahin, nsahin, vu}.

It is easy to observe that a monomial curve in 2-space has the same resolution with its tangent cone and the Betti sequence is $(1,1)$. A monomial curve in 3-space has the Betti sequence $(1,2,1)$, if it is a complete intersection. Otherwise, it has the Betti sequence $(1,3,2)$ by a result of Herzog \cite{paper17}, \cite{stamate}. Moreover, if it has a Cohen-Macaulay tangent cone, then the Betti sequence of the tangent cone is again $(1,2,1)$ or $(1,3,2)$. This is a consequence of a result given both in \cite{herzog0} and \cite{robbiano}, stating that a monomial curve $C$ in 3-space has a Cohen-Macaulay tangent cone if and only if the minimal number of generators of the defining ideal of the tangent cone is less than or equal to 3.  The generators of the ideal of the tangent cone of the monomial curve in 3-space can be arbitrarily large and Shibuta gives a family with arbitrary large minimal generators \cite{goto}. Herzog and Stamate extend the example of Shibuta to obtain families of monomial curves having tangent cones with arbitrary large Betti numbers and compute them \cite{herzog1}.  

In 4-space, there are monomial curve families with arbitrary large minimal generators given by authors     
\cite{paper01}, \cite{bresinsky}. Stamate computes the Betti numbers for these families and their tangent cones in \cite{stamate}.  
If a monomial curve in $4$-space is a complete intersection, then the Betti sequence is $(1,3,3,1)$. If it is not a complete intersection monomial curve with the corresponding numerical semigroup symmetric, then Barucci, Fr{\"o}berg and Şahin showed that the Betti sequence of the monomial curve is $(1,5,5,1)$ \cite{barucci}. Bresinsky showed that a non-complete intersection monomial curve in 4-space with the corresponding numerical semigroup symmetric has always 5 generators and explicitly described the generator set \cite{paper07}. By using this description, Mete and Zengin \cite{paper21} and Mete \cite{paper22} showed that possible Betti sequences for a Cohen-Macaulay tangent cone are $(1,5,6,2), (1,5,5,1)$ and $(1,6,8,3)$.  

\section{Nice gluing}
 
We first give the definitions of a gluing and a nice gluing. The concept of nice gluing introduced in \cite{Arslan1} is an effective tool to form families of complete intersection monomial curves with complete intersection tangent cones.

\begin{definition} \cite{rosales}
Let $S_1=\langle n_1, n_2, ..., n_k \rangle$ and $S_2=\langle m_1, m_2, ..., m_l \rangle$ be two numerical semigroups minimally generated by $n_1 < n_2 < \cdots < n_k $ and $m_1 < m_2 < \cdots < m_l $ respectively.
Let
$p=a_1n_1+\dots+a_kn_k \in S_1$ and $q=b_1m_1+\dots+b_lm_l \in S_2$
be two positive integers satisfying $gcd(p,q)=1$ with $p \not\in
\{m_{1},\dots,m_{l}\}$ and $q \not\in \{n_{1},\dots,n_{k}\}$ and $\{pm_{1},\dots,pm_{l}\} \cap \{qn_{1},\dots,qn_{k}\}=\emptyset$. The
numerical semigroup $S=\langle
qn_{1},\dots,qn_{k},pm_{1},\dots,pm_{l} \rangle$ is called a gluing
of the semigroups $S_1$ and $S_2$.
\end{definition} 

\begin{definition} \label{nice} \cite{Arslan1}
A gluing $S=\langle
qn_{1},\dots,qn_{k},pm_{1},\dots,pm_{l} \rangle$ 
of the semigroups $S_1$ and $S_2$ is called a nice gluing, if 
$p =a_1 n_1 + a_2 n_2 + ...+ a_k n_k$ and $q = b_1 m_1 $ with $b_1 \le a_1 + a_2 + ...+ a_k $. 
\end{definition} 

Delorme's characterization reinterpreted by Rosales in terms of gluing says that every complete intersection numerical semigroup is a gluing of two complete intersection numerical semigroups \cite{delorme}, \cite{rosales}.
We first prove that having a complete intersection tangent cone for a monomial curve corresponding to a  
numerical semigroup $\langle c_1, c_2, c_3 \rangle$ is equivalent to being a numerical semigroup obtained from a nice gluing.    
\begin{theorem}
The numerical symmetric semigroup $S= \langle c_1, c_2, c_3 \rangle$ is obtained by a nice gluing if and only if the monomial curve $C$ with parametrization $x_1=t^{c_1},x_2=t^{c_2},x_3=t^{c_3}$ has a complete intersection tangent cone.
\end{theorem}

\begin{proof}
First, assume that $S= \langle c_1, c_2, c_3 \rangle$ is obtained by a nice gluing with $c_1$ less than $c_2$ and $c_3$. Without loss of generality, we can assume that either $gcd(c_1,c_2)$ or $gcd(c_2,c_3)$ is greater than 1, since if $gcd(c_1,c_3)$ is greater than 1, we can simply interchange $c_2$ and $c_3$.  
Hence, $S= \langle c_1, c_2, c_3 \rangle$ is either a nice gluing of $\langle n_1, n_2 \rangle$ with $n_1<n_2$ and $\langle 1 \rangle$ or a nice gluing of $\langle 1 \rangle$ and $\langle m_1, m_2 \rangle$ with $m_1<m_2$.
 \begin{enumerate}
     \item By using the notation in Definition \ref{nice}, $\langle c_1, c_2, c_3 \rangle$ is a nice gluing of $\langle n_1, n_2 \rangle$ and $\langle m_1 \rangle=\langle 1 \rangle$ means $c_1=qn_1, \; c_2=qn_2, \; c_3=pm_1=p$ with $p=a_1n_1+a_2n_2$, $q=b_1m_1=b_1$ and $b_1 \leq a_1+a_2$. Then the defining ideal $I(C)$ is generated by $\{x_2^{n_1} - x_1^{n_2},x_3^{b_1} - x_1^{a_1}x_2^{a_2}\}$. The leading monomials are $LM(x_1^{n_2} - x_2^{n_1})=x_2^{n_1}$ and $LM(x_3^{b_1} - x_1^{a_1}x_2^{a_2})=x_3^{b_1}$ with respect to the negative degree reverse
lexicographical ordering with $x_3>x_2>x_1$. Since the leading monomials are relatively prime, this is a standard basis. From \cite[Lemma 5.5.11]{greuel}, if we have a standard basis for $I(C)$ with respect to a local ordering, then the tangent cone $I(C)_*$ is generated by the least homogeneous summands of the elements in the standard basis. Hence, in our case, the tangent cone is generated by two elements and it is a complete intersection. 
    \item $\langle c_1, c_2, c_3 \rangle$ is a nice gluing of $\langle n_1 \rangle = \langle 1 \rangle$ and $\langle m_1, m_2 \rangle$ means $c_1=qn_1=q, \; c_2=pm_1, \; c_3=pm_2$ with $p=a_1n_1=a_1$, $q=b_1m_1$ and $b_1 < a_1$. Then the defining ideal $I(C)$ is generated by $\{x_2^{m_2} - x_3^{m_1},x_1^{a_1} - x_2^{b_1}\}$. The leading monomials are $LM(x_2^{m_2} - x_3^{m_1})=x_3^{m_1}$ and $LM(x_1^{a_1} - x_2^{b_1})=x_2^{b_1}$ with respect to any local ordering. Since the leading monomials are relatively prime, this is a standard basis. With the same argument above, the tangent cone $I(C)_*$ is generated by these two monomials, so it is also a complete intersection.
\end{enumerate}
Conversely, let the symmetric monomial curve $C$ defined by $x_1=t^{c_1},x_2=t^{c_2},x_3=t^{c_3}$ with $c_1 < c_2 < c_3 $ have a complete intersection tangent cone. It was shown in \cite[Corollary 3.2]{Shen} by Shen that the tangent cone is a complete intersection if either $I(C)=\langle x_1^{\alpha_1}-x_2^{\alpha_2},x_3^{\alpha_3}-x_1^{\alpha_{31}}x_2^{\alpha_{32}} \rangle$ or
$I(C)=\langle x_1^{\alpha_1}-x_3^{\alpha_3},x_2^{\alpha_2}-x_1^{\alpha_{21}}x_3^{\alpha_{23}} \rangle$ with $\alpha_2 \leq \alpha_{21}+\alpha_{23}$. Also note that, $c_1 < c_2 < c_3$ implies $\alpha_2 < \alpha_1$, $\alpha_3 < \alpha_1$ and $\alpha_3 < \alpha_{31}+\alpha_{32}$.  
 \begin{enumerate}
     \item If $I(C)=\langle x_1^{\alpha_1}-x_2^{\alpha_2},x_3^{\alpha_3}-x_1^{\alpha_{31}}x_2^{\alpha_{32}} \rangle$, then $\langle c_1, c_2, c_3 \rangle$ is a gluing of $\langle \alpha_2, \alpha_1 \rangle$ and 
$\langle 1 \rangle$ with $c_1=\alpha_3 \alpha_2$, $c_2=\alpha_3 \alpha_1$ and $c_3=\alpha_{31}\alpha_{2}+\alpha_{32}\alpha_{1}$. Since $\alpha_3 < \alpha_{31}+\alpha_{32}$, this is a nice gluing.  
\item If $I(C)=\langle x_1^{\alpha_1}-x_3^{\alpha_3},x_2^{\alpha_2}-x_1^{\alpha_{21}}x_3^{\alpha_{23}} \rangle$ with $\alpha_2 \leq \alpha_{21}+\alpha_{23}$, then $\langle c_1, c_2, c_3 \rangle$ is a gluing of $\langle \alpha_3, \alpha_1 \rangle$ and 
$\langle 1 \rangle$ with $c_1=\alpha_2 \alpha_3$, $c_3=\alpha_2 \alpha_1$ and $c_2=\alpha_{21}\alpha_{3}+\alpha_{23}\alpha_{1}$. The condition $\alpha_2 \leq \alpha_{21}+\alpha_{23}$ implies that this is a nice gluing.  
\end{enumerate}
\end{proof}

\begin{remark}
Note that the complete intersection numerical semigroup $\langle 6, 7, 15 \rangle$ is obtained by gluing, which is not nice, and the tangent cone of the monomial curve is not a complete intersection. Moreover, it is not even Cohen-Macaulay.   
\end{remark}

\begin{remark}
Theorem 2.3 can also be interpreted within the framework of  
the $\gamma$-rectangular Apéry sets introduced by D’Anna, Micale and Sammartano in \cite{D'Anna} and elaborated in \cite{D'Anna2}. To this end, we will first review some necessary definitions. For a numerical semigroup $S$ and for any $n \in S$, the set
\[
\operatorname{Ap}(S,n) = \{\, s \in S \mid s-n \notin S \,\}
\]
is called the Apéry set of $S$ with respect to $n$. Equivalently, it is the set of the smallest elements of $S$ in each congruence class modulo $n$. The smallest positive integer in $S$ is called the multiplicity of $S$, denoted by $m(S)$.  
When $n = m(S)$, we simply write the Apéry set as $\operatorname{Ap}(S)$. A representation of $s \in S$ is an expression of the form
\[
s = a_1 n_1 + \cdots + a_l n_l, \quad a_i \in \mathbb{N}.
\]
Given $s \in S$, we define the \emph{order} of $s$ by
\[
\operatorname{ord}(s) = \max \Big\{ a_1 + \cdots + a_l \;\Big|\; 
s = a_1 n_1 + \cdots + a_l n_l \text{ is a representation of } s \Big\}.
\]
A representation is called \emph{maximal} if 
$a_1 + \cdots + a_l = \operatorname{ord}(s)$
and it is called a \emph{unique maximal representation}, if there is only one way 
to write $s$ as a linear combination of $n_1, \dots, n_l$ achieving this maximum. In \cite{D'Anna},
a numerical semigroup $S=\langle n_1, \dots, n_l \rangle$ with multiplicity $n_1$ is said to have an $\alpha$-rectangular Apéry set if
\[
\operatorname{Ap}(S) = 
\Big\{ \sum_{i=2}^{l} a_i n_i \ \Big|\ 
0 \leq a_i \leq \alpha_i \Big\},
\]
where
\begin{align*}
\alpha_i = \max \{\, h \in \mathbb{N} \;&\;\Big|\;
h n_i \in \operatorname{Ap}(S)\},
\end{align*}
and
$\gamma$-rectangular Apéry set if
\[
\operatorname{Ap}(S) = 
\Big\{ \sum_{i=2}^{l} a_i n_i \ \Big|\ 
0 \leq a_i \leq \gamma_i \Big\},
\]
where
\begin{align*}
\gamma_i = \max \{\, h \in \mathbb{N} \;&\;\Big|\;
h n_i \in \operatorname{Ap}(S),\;
\operatorname{ord}(h n_i)=h, \\
&\;\text{and $h n_i$ has a unique maximal representation} \,\}.
\end{align*}
In \cite{D'Anna2}, D'Anna et al. prove the following theorem:

\begin{theorem} \cite[Theorem 3.3]{D'Anna2}
Let $S_1 = \langle n_1, n_2,..., n_l \rangle$ be a numerical semigroup 
having an $\alpha$-rectangular Apéry set. 
Let $S_2 = \mathbb{N}$. Then, the gluing $S$ of the semigroups $S_1$ and $S_2$ 
given by $S=\langle
qn_{1},\dots,qn_{l},p \rangle$ with $p=a_1n_1+...+a_ln_l \notin \operatorname{Ap}(S_1)$ and $qn_1 < p$ has an $\alpha$-rectangular Apéry set. Conversely, every semigroup $S \neq \mathbb{N}$ with an $\alpha$-rectangular Apéry set arises in this way.
\end{theorem}

Note that the condition $qn_1 < p$ makes this gluing a nice gluing, but it has an extra condition $p \notin \operatorname{Ap}(S_1)$. Since every $\alpha$-rectangular Apéry set is also a $\gamma$-rectangular Apéry set, the theorem is also true for semigroups having $\gamma$-rectangular Apéry sets. It is natural to ask if we can drop the extra condition $p \notin \operatorname{Ap}(S_1)$ for semigroups with $\gamma$-rectangular Apéry set:

\begin{question}
Let $S_1 = \langle n_1, n_2,..., n_l \rangle$ be a numerical semigroup 
having a $\gamma$-rectangular Apéry set. 
Let $S_2 = \mathbb{N}$. Does the nice gluing $S$ of the semigroups $S_1$ and $S_2$ 
given by $S=\langle
qn_{1},\dots,qn_{l},p \rangle$ with $p=a_1n_1+...+a_ln_l$ and $q \leq a_1+...+a_l$ have a $\gamma$-rectangular Apéry set? 
\end{question}  

\begin{example}
First, note that every numerical semigroup $\langle n_1, n_2 \rangle$ has an $\alpha$-rectangular Apéry set so that it also has a $\gamma$-rectangular Apéry set. The nice gluing of $S_1=\langle 4, 5 \rangle$ and $S_2 = \mathbb{N}$ with $q=2$ and $p=15$ gives the numerical semigroup $S=\langle 8, 10, 15 \rangle$. Here, $p \in \operatorname{Ap}(S_1)$ and $S$ does not have an $\alpha$-rectangular Apéry set, but it has a $\gamma$-rectangular Apéry set.      
\end{example}

Theorem 2.3 can be considered as a partial answer to Question 2.7 and one direction of the theorem can be stated as follows:

\begin{theorem}
The nice gluing $S$ of the semigroup $S_1 = \langle n_1, n_2 \rangle$ and $S_2 = \mathbb{N}$ 
given by $S=\langle
qn_{1},qn_{2},p \rangle$ with $p=a_1n_1+a_2n_2$ and $q \leq a_1+a_2$ has a $\gamma$-rectangular Apéry set. 
\end{theorem}

\begin{proof}
Theorem 2.3 implies that $S$ has a complete intersection tangent cone, and it follows from \cite[Theorem 3.6]{D'Anna} that $S$ has a $\gamma$-rectangular Apéry set. 
\end{proof}

\end{remark}

We can use the idea in Theorem 2.3 to obtain infinitely many families of complete intersection monomial curves with complete intersection tangent cones in 4-space so that the Betti sequence is $(1, 3, 3, 1)$ for all of them.  

\begin{theorem}
 If the symmetric numerical semigroup $S= \langle c_1, c_2, c_3, c_4 \rangle$ is obtained by a nice gluing of 
\begin{enumerate}  
 \item $\langle n_1, n_2 \rangle$ and $\langle m_1, m_2 \rangle$, or
\item symmetric $\langle n_1, n_2, n_3 \rangle$ and $\langle 1 \rangle$ (where the complete intersection monomial curve defined by $x_1=t^{n_1},x_2=t^{n_2},x_3=t^{n_3}$ has a complete intersection tangent cone), or
\item $\langle 1 \rangle$ and symmetric $\langle m_1, m_2, m_3 \rangle$ (where the complete intersection monomial curve defined by  $x_1=t^{m_1},x_2=t^{m_2},x_3=t^{m_3}$ has a complete intersection tangent cone)
\end{enumerate}
then the monomial curve $C$ defined by 
$x_1=t^{c_1},x_2=t^{c_2},x_3=t^{c_3},x_4=t^{c_4}$ has a complete intersection tangent cone.
\end{theorem}

\begin{proof}
\begin{enumerate} 
\item $\langle c_1, c_2, c_3, c_4 \rangle$ is nice gluing of $\langle n_1, n_2 \rangle$ and $\langle m_1, m_2 \rangle$ means $c_1=b_1m_1n_1$, $c_2=b_1m_1n_2, c_3=(a_1n_1+a_2n_2)m_1$ and $c_4=(a_1n_1+a_2n_2)m_2$ where $b_1 \leq a_1+a_2$. Then the defining ideal $I(C)$ is generated by $\{x_2^{n_1} - x_1^{n_2},x_4^{m_1} - x_3^{m_2},x_3^{b_1} - x_1^{a_1}x_2^{a_2}\}$. We have $LM(x_2^{n_1} - x_1^{n_2})=x_2^{n_1}$, $LM(x_4^{m_1} - x_3^{m_2})=x_4^{m_1}$ and $LM(x_3^{b_1} - x_1^{a_1}x_2^{a_2})=x_3^{b_1}$ with respect to the negative degree reverse
lexicographical ordering with $x_4>x_3>x_2>x_1$. Since any two of the monomials are relatively prime, this set is a standard basis. Hence, from \cite[Lemma 5.5.11]{greuel}, the tangent cone is also generated by three elements, so it is also a complete intersection.
\item If the complete intersection monomial curve given by $x_1=t^{n_1},x_2=t^{n_2},x_3=t^{n_3}$ with $n_1<n_2<n_3$ has a complete intersection tangent cone, we know from Shen's classification that the defining ideal is either $\langle x_1^{\alpha_1}-x_2^{\alpha_2},x_3^{\alpha_3}-x_1^{\alpha_{31}}x_2^{\alpha_{32}} \rangle$ or
$\langle x_1^{\alpha_1}-x_3^{\alpha_3},x_2^{\alpha_2}-x_1^{\alpha_{21}}x_3^{\alpha_{23}} \rangle$ with $\alpha_2 \leq \alpha_{21}+\alpha_{23}$, \cite[Corollary 3.2]{Shen}. A nice gluing of $\langle n_1, n_2, n_3 \rangle$ and $\langle 1 \rangle$ corresponds to the monomial curve $C$ with the generator sets $S_1=\{ x_1^{\alpha_1}-x_2^{\alpha_2},x_3^{\alpha_3}-x_1^{\alpha_{31}}x_2^{\alpha_{32}},x_4^{b_1}-x_1^{a_1}x_2^{a_2}x_3^{a_3} \}$ in the former case and $S_2=\{ x_1^{\alpha_1}-x_3^{\alpha_3},x_2^{\alpha_2}-x_1^{\alpha_{21}}x_3^{\alpha_{23}}, x_4^{b_1}-x_1^{a_1}x_2^{a_2}x_3^{a_3} \}$ in the latter case, where  
$b_1 \leq a_1+a_2+a_3$. In both cases, $S_1$ and $S_2$ are standard bases with respect to the negative degree reverse
lexicographical ordering with $x_4>x_2>x_3>x_1$. Hence, the tangent cone is generated by 3 elements and it is a complete intersection.
   
\item The proof can be done as in (2).
\end{enumerate}

\end{proof} 

We have shown that a complete intersection monomial curve corresponding to a numerical semigroup obtained by a nice gluing of two semigroups with corresponding monomial curves having complete intersection tangent cones has also a complete intersection tangent cone. It is natural to ask if the converse holds:

\begin{question} If the complete intersection monomial curve $C$ defined by 
$x_1=t^{c_1},x_2=t^{c_2},x_3=t^{c_3},x_4=t^{c_4}$ has a complete intersection tangent cone, then is $\langle c_1, c_2, c_3, c_4 \rangle$ a nice gluing of two numerical semigroups with corresponding monomial curves having complete intersection tangent cones?
\end{question} 

It is also natural to ask the same question for monomial curves in arbitrary dimensions:  

\begin{question} Does the complete intersection monomial curve $C$ defined by 
$x_1=t^{c_1},x_2=t^{c_2},...,x_n=t^{c_n}$ have a complete intersection tangent cone if and only if $\langle c_1, c_2,..., c_n \rangle$ is a nice gluing of two numerical semigroups with corresponding monomial curves having complete intersection tangent cones?
\end{question} 

In \cite[Theorem 2.14]{herzog2}, Herzog and Stamate give a family of numerical semigroups obtained by quadratic gluings such that $H_0=\mathbb{N}, H_1=\langle 2H_0, l_1 \rangle,...,H_{n}=\langle 2H_{n-1}, l_{n} \rangle$, where $l_i$ is an odd integer greater than 1. In fact, this corresponds to a sequence of nice gluings. The defining ideal of the complete intersection monomial curve corresponding to the numerical semigroup $H_n$ is generated by the set $$\{x_2^2-x_1^{l_{11}},x_3^2-x_1^{l_{21}}x_2^{l_{22}},...,x_{n+1}^2-x_1^{l_{n1}}x_2^{l_{n2}}...x_n^{l_{nn}}\}$$
where $l_i=l_{i1}+l_{i2}+...+l_{ii}$ for $1 \leq i \leq n$ and $l_i$ is an odd integer greater than 1. It is straightforward to check that this set is a standard basis with respect to the negative degree reverse
lexicographical ordering with $x_{n+1}>x_n>...>x_1$, which shows that it has a complete intersection tangent cone. Thus, this family is a nice example of complete intersection monomial curves obtained by nice gluings and having complete intersection tangent cones.     

\section{Families obtained by non-nice gluing}

In this section, we give explicit minimal free resolutions of the tangent cones of a family of complete intersection monomial curves in an affine 4-space obtained by a non-nice gluing.  

We use the families constructed in \cite{Arslan2}. By a nice gluing of $S=\langle n_1, n_2 \rangle$ and $\langle 1 \rangle$, they first obtain the semigroup $S'=\langle pn_1, pn_2, a_1n_1+a_2n_2 \rangle$ with $p \leq a_1+a_2$ and then by a non-nice gluing of $\langle 1 \rangle$ and $S'$, they construct a complete intersection numerical semigroup family 
\begin{center}
    $ \langle (s-d)pn_1 + pn_2, spn_1, spn_2, s(a_1n_1+a_2n_2) \rangle$
\end{center}
with $p\leq a_1 + a_2$, $gcd(p,a_1n_1 + a_2n_2)=1$, $gcd((s-d)pn_1 + pn_2, s)=1$, $s>d$, $d= \lceil \frac{n_2}{n_1} \rceil$ and $(s-d)pn_1 + pn_2$ is the smallest among the generators, \cite{Arslan2}. The numerical semigroup
corresponds to a complete intersection monomial curve with the parametrization $y_1= t^{(s-d)pn_1 + pn_2}, x_1=t^{spn_1}, x_2=t^{spn_2}, x_3=t^{s(a_1n_1+a_2n_2)}$ and they show that the minimal standard basis of this monomial curve is
\begin{center}
    $ \{ f_1=x_1^{n_2} - x_2^{n_1}, f_2=x_3^p - x_1^{a_1}x_2^{a_2}, f_3=y_1^s - x_1^{s-d}x_2, f_{j+3}=y_1^{js}x_2^{n_1 - j} - x_1^{n_2 +j(s-d)}\}$
\end{center}
for $1\leq j \leq i_S= \lceil \frac{n_2-n_1}{d-1} \rceil$ with respect to the negative lexicographical ordering with $x_3>x_2>x_1>y_1$ \cite{Arslan2}.

Moreover, by using this construction, they start with $\langle n_1, n_2 \rangle = \langle 8, 27 \rangle$ which leads to $d=4$ and choose $a_1=p$ and $a_2=1$ to construct a family $ \langle (s-4)8p + 27p, 8sp, 27sp, (8p+27)s \rangle$ with two parameters \cite{Arslan2}.
In this section, we give the explicit minimal free resolutions of the tangent cones of the monomial curves $C_{s,p}$ corresponding to this numerical semigroup family. The computations are done in \cite{macaulay2} and \cite{singular}. 
 
\begin{theorem}

Let $C_{s,p} = C((s-4)8p+27p, 8sp, 27sp, (8p+27)s)$ be a complete intersection monomial curve, where $gcd(p,3)=1$ and $gcd(8p(s-4)+27p, s)=1$. Then the defining ideal of the tangent cone is
$$I(C_{s,p})_* = (x_2^8, x_3^p, x_1^{s-4}x_2, y_1^sx_2^7, y_1^{2s}x_2^6, y_1^{3s}x_2^5, y_1^{4s}x_2^4, y_1^{5s}x_2^3, y_1^{6s}x_2^2, x_1^{7s-1}),$$
and the sequence of $R$-modules, where $R = K[y_1,x_1,x_2,x_3]$,
    $$  
    0 \longrightarrow R^6 \xrightarrow{\varphi_4} R^{20} \xrightarrow{\varphi_3} R^{23} \xrightarrow{\varphi_2} R^{10} \xrightarrow{\varphi_1} R^1 \rightarrow R/I(C_{s,p})_* \longrightarrow 0
    $$
is a minimal free resolution for the tangent cone of $C_{s,p}$, where
    $$
    \varphi_1=\begin{pmatrix} x_3^p & x_1^{s-4}x_2 &x_2^8 & y_1^sx_2^7 & y_1^{2s}x_2^6 & y_1^{3s}x_2^5 & y_1^{4s}x_2^4 & y_1^{5s}x_2^3 & y_1^{6s}x_2^2 & x_1^{7s-1}\end{pmatrix}
    $$

\begin{landscape}
\thispagestyle{empty}
   $\varphi_2=\begin{psmallmatrix} -x_1^{s-4}x_2 & -x_2^8 & 0 & 0 &-x_2^7y_1^s & 0 & 0 & -x_2^6y_1^{2s} & 0 & 0 & -x_2^5y_1^{3s} & 0 & 0 & -x_2^4y_1^{4s} & 0 & 0 & -x_2^3y_1^{5s} & 0 & 0 & -x_2^2y_1^{6s} & 0 & 0 & -x_1^{7s-1} \\ 
    x_3^p & 0 & -x_2^7 & 0 & 0 & -x_2^6y_1^{s} & 0 & 0 & -x_2^5y_1^{2s} & 0 & 0 & -x_2^4y_1^{3s} & 0 & 0 & -x_2^3y_1^{4s} & 0 & 0 & -x_2^2y_1^{5s} & 0 & 0 & -x_2y_1^{6s} & -x_1^{6s+3} & 0 \\ 
    0 & x_3^p & x_1^{s-4} & -y_1^s & 0 & 0 & 0 & 0 & 0 & 0 & 0 & 0 & 0 & 0 & 0 & 0 & 0 & 0 & 0 & 0 & 0 & 0 & 0 \\ 
    0 & 0 & 0 & x_2 & x_3^p & x_1^{s-4} & -y_1^s & 0 & 0 & 0 & 0 & 0 & 0 & 0 & 0 & 0 & 0 & 0 & 0 & 0 & 0 & 0 & 0 \\ 
    0 & 0 & 0 & 0 & 0 & 0 & x_2 & x_3^p & x_1^{s-4} & -y_1^s & 0 & 0 & 0 & 0 & 0 & 0 & 0 & 0 & 0 & 0 & 0 & 0 & 0 \\ 
    0 & 0 & 0 & 0 & 0 & 0 & 0 & 0 & 0 & x_2 & x_3^p & x_1^{s-4} & -y_1^s & 0 & 0 & 0 & 0 & 0 & 0 & 0 & 0 & 0 & 0 \\ 
    0 & 0 & 0 & 0 & 0 & 0 & 0 & 0 & 0 & 0 & 0 & 0 & x_2 & x_3^p & x_1^{s-4} & -y_1^s & 0 & 0 & 0 & 0 & 0 & 0 & 0 \\ 
    0 & 0 & 0 & 0 & 0 & 0 & 0 & 0 & 0 & 0 & 0 & 0 & 0 & 0 & 0 & x_2 & x_3^p & x_1^{s-4} &-y_1^s & 0 & 0 & 0 & 0 \\ 
    0 & 0 & 0 & 0 & 0 & 0 & 0 & 0 & 0 & 0 & 0 & 0 & 0 & 0 & 0 & 0 & 0 & 0 & x_2 & x_3^p & x_1^{s-4} & 0 & 0 \\ 
    0 & 0 & 0 & 0 & 0 & 0 & 0 & 0 & 0 & 0 & 0 & 0 & 0 & 0 & 0 & 0 & 0 & 0 & 0 & 0 & 0 & x_2 & x_3^p \end{psmallmatrix}$
    
    $\varphi_3= \begin{psmallmatrix} x_2^7 & 0 & 0 & x_2^6y_1^s & 0 & 0 & x_2^5y_1^{2s} & 0 & 0 & x_2^4y_1^{3s} & 0 & 0 & x_2^3y_1^{4s} & 0 & 0 & x_2^2y_1^{5s} & 0 & 0 & x_2y_1^{6s} & -x_1^{6s+3} \\ 
    -x_1^{s-4} & -y_1^s & 0 & 0 & 0 & 0 & 0 & 0 & 0 & 0 & 0 & 0 & 0 & 0 & 0 & 0 & 0 & 0 & 0 & 0 \\  
    x_3^p & 0 & -y_1^s & 0 & 0 & 0 & 0 & 0 & 0 & 0 & 0 & 0 & 0 & 0 & 0 & 0 & 0 & 0 & 0 & 0 \\  
    0 & -x_3^p & -x_1^{s-4} & 0 & 0 & 0 & 0 & 0 & 0 & 0 & 0 & 0 & 0 & 0 & 0 & 0 & 0 & 0 & 0 & 0 \\  
    0 & x_2 & 0 & -x_1^{s-4} & -y_1^s & 0 & 0 & 0 & 0 & 0 & 0 & 0 & 0 & 0 & 0 & 0 & 0 & 0 & 0 & 0 \\  
    0 & 0 & x_2 & x_3^p & 0 & -y_1^s & 0 & 0 & 0 & 0 & 0 & 0 & 0 & 0 & 0 & 0 & 0 & 0 & 0 & 0\\  
    0 & 0 & 0 & 0 & -x_3^p & -x_1^{s-4} & 0 & 0 & 0 & 0 & 0 & 0 & 0 & 0 & 0 & 0 & 0 & 0 & 0 & 0 \\  
    0 & 0 & 0 & 0 & x_2 & 0 & -x_1^{s-4} & -y_1^s & 0 & 0 & 0 & 0 & 0 & 0 & 0 & 0 & 0 & 0 & 0 & 0 \\  
    0 & 0 & 0 & 0 & 0 & x_2 & x_3^p & 0 & -y_1^s & 0 & 0 & 0 & 0 & 0 & 0 & 0 & 0 & 0 & 0 & 0 \\  
    0 & 0 & 0 & 0 & 0 & 0 & 0 & -x_3^p & -x_1^{s-4} & 0 & 0 & 0 & 0 & 0 & 0 & 0 & 0 & 0 & 0 & 0 \\  
    0 & 0 & 0 & 0 & 0 & 0 & 0 & x_2 & 0 & -x_1^{s-4} & -y_1^s & 0 & 0 & 0 & 0 & 0 & 0 & 0 & 0 & 0 \\  
    0 & 0 & 0 & 0 & 0 & 0 & 0 & 0 & x_2 & x_3^p & 0 & -y_1^s & 0 & 0 & 0 & 0 & 0 & 0 & 0 & 0 \\  
    0 & 0 & 0 & 0 & 0 & 0 & 0 & 0 & 0 & 0 & -x_3^p & -x_1^{s-4} & 0 & 0 & 0 & 0 & 0 & 0 & 0 & 0 \\
    0 & 0 & 0 & 0 & 0 & 0 & 0 & 0 & 0 & 0 & x_2 & 0 & -x_1^{s-4} & -y_1^s & 0 & 0 & 0 & 0 & 0 & 0 \\  
    0 & 0 & 0 & 0 & 0 & 0 & 0 & 0 & 0 & 0 & 0 & x_2 & x_3^p & 0 & -y_1^s & 0 & 0 & 0 & 0 & 0 \\  
    0 & 0 & 0 & 0 & 0 & 0 & 0 & 0 & 0 & 0 & 0 & 0 & 0 & -x_3^p & -x_1^{s-4} & 0 & 0 & 0 & 0 & 0 \\  
    0 & 0 & 0 & 0 & 0 & 0 & 0 & 0 & 0 & 0 & 0 & 0 & 0 & x_2 & 0 & -x_1^{s-4} & -y_1^s & 0 & 0 & 0 \\  
    0 & 0 & 0 & 0 & 0 & 0 & 0 & 0 & 0 & 0 & 0 & 0 & 0 & 0 & x_2 & x_3^p & 0 & -y_1^s & 0 & 0 \\  
    0 & 0 & 0 & 0 & 0 & 0 & 0 & 0 & 0 & 0 & 0 & 0 & 0 & 0 & 0 & 0 & -x_3^p & -x_1^{s-4} & 0 & 0 \\  
    0 & 0 & 0 & 0 & 0 & 0 & 0 & 0 & 0 & 0 & 0 & 0 & 0 & 0 & 0 & 0 & x_2 & 0 & -x_1^{s-4} & 0 \\  
    0 & 0 & 0 & 0 & 0 & 0 & 0 & 0 & 0 & 0 & 0 & 0 & 0 & 0 & 0 & 0 & 0 & x_2 & x_3^p & 0 \\ 
    0 & 0 & 0 & 0 & 0 & 0 & 0 & 0 & 0 & 0 & 0 & 0 & 0 & 0 & 0 & 0 & 0 & 0 & 0 & -x_3^p \\ 
    0 & 0 & 0 & 0 & 0 & 0 & 0 & 0 & 0 & 0 & 0 & 0 & 0 & 0 & 0 & 0 & 0 & 0 & 0 & x_2 \\ \end{psmallmatrix}$
    
\end{landscape}

    $$
    \varphi_4=\begin{pmatrix} -y_1^s & 0 & 0 & 0 & 0 & 0 \\
    x_1^{s-4} & 0 & 0 & 0 & 0 & 0 \\
    -x_3^p & 0 & 0 & 0 & 0 & 0 \\
    x_2 & -y_1^s & 0 & 0 & 0 & 0 \\
    0 & x_1^{s-4} & 0 & 0 & 0 & 0 \\
    0 & -x_3^p & 0 & 0 & 0 & 0 \\
    0 & x_2 & -y_1^s & 0 & 0 & 0 \\
    0 & 0 & x_1^{s-4} & 0 & 0 & 0 \\
    0 & 0 & -x_3^p & 0 & 0 & 0 \\
    0 & 0 & x_2 & -y_1^s & 0 & 0 \\
    0 & 0 & 0 & x_1^{s-4} & 0 & 0 \\
    0 & 0 & 0 & -x_3^p & 0 & 0 \\
    0 & 0 & 0 & x_2 & -y_1^s & 0 \\
    0 & 0 & 0 & 0 & x_1^{s-4} & 0 \\
    0 & 0 & 0 & 0 & -x_3^p & 0 \\
    0 & 0 & 0 & 0 & x_2 & -y_1^s \\
    0 & 0 & 0 & 0 & 0 & x_1^{s-4} \\
    0 & 0 & 0 & 0 & 0 & -x_3^p \\
    0 & 0 & 0 & 0 & 0 & x_2 \\
    0 & 0 & 0 & 0 & 0 & 0 \\
    \end{pmatrix}
    $$

\end{theorem}

\begin{proof}
We use the technique used in \cite{barucci}, \cite{paper21} and \cite{paper22}. It is a direct computation to check that $\varphi_1 \varphi_2=\varphi_2 \varphi_3=\varphi_3 \varphi_4 = 0$. Thus, the sequence is a complex. By using Buchsbaum-Eisenbud criterion \cite[Corollary 2]{Buch}, we will show the exactness. We need to check that $rank(\varphi_1) + rank(\varphi_2)=10$, $rank(\varphi_2) + rank(\varphi_3) = 23$ and $rank(\varphi_3) + rank(\varphi_4) = 20$ and $I(\varphi_i)$ which is the ideal generated by the minors of maximal rank of $\varphi_i$ contains an $R$-sequence of length $i$.

It is clear that $rank(\varphi_1)=1$. Every $10 \times 10$ minor of $\varphi_2$ is 0. In the matrix $\varphi_2$, if we delete the 1st row and the 3rd, 5th, 6th, 8th, 9th, 11th, 12th, 14th, 15th, 17th, 18th, 20th, 21st and 23rd columns, we have $x_2^7x_3^{2p}$, and if we delete the 10th row and the 1st, 2nd, 3rd, 5th, 6th, 8th, 9th, 11th, 12th, 14th, 15th, 17th, 18th and 20th columns, we have $-x_1^{14s-2}y_1^{6s}$. These two $9 \times 9$ minors are relatively prime, so $I(\varphi_2)$ has a regular sequence of length 2. 

Every $15 \times 15$ minor of $\varphi_3$ is 0. In the matrix $\varphi_3$, if we delete the 4th, 7th, 10th, 13th, 16th, 19th, 20th, 21st and 22nd rows and the 1st, 4th, 7th, 10th, 13th and 16th columns, we have $x_2^2y_1^{18s}$, if we delete the 1st, 2nd, 5th, 8th, 11th, 14th, 17th, 20th and 23rd rows and the 3rd, 6th, 9th, 12th, 15th and 18th columns, we have $-x_3^{14p}$, and if we delete the 3rd, 6th, 9th, 12th, 15th, 18th, 21st, 22nd and 23rd rows and the 2nd, 5th, 8th, 11th, 14th and 17th columns, we have $-x_1^{19s-49}$. These three $14 \times 14$ minors are relatively prime, so $I(\varphi_3)$ has a regular sequence of length 3.

In the matrix $\varphi_4$, if we delete the 2nd, 3rd, 5th, 6th, 8th, 9th, 11th, 12th, 14th, 15th, 17th, 18th, 19th and 20th rows, we have $y_1^{6s}$, if we delete the 1st, 3rd, 4th, 6th, 7th, 9th, 10th, 12th, 13th, 15th, 16th, 18th, 19th and 20th rows, we have $x_1^{6s-24}$, if we delete the 1st, 2nd, 4th, 5th, 7th, 8th, 10th, 11th, 13th, 14th, 16th, 17th, 19th and 20th rows, we have $x_3^{6p}$ and if we delete the 1st, 2nd, 3rd, 5th, 6th, 8th, 9th, 11th, 12th, 14th, 15th, 17th, 18th and 20th rows, we have $x_2^6$. These four $6 \times 6$ minors are relatively prime, so $I(\varphi_4)$ has a regular sequence of length 4.
\end{proof}

Hence, we have constructed families of complete intersection monomial curves with two free parameters and computed the minimal free resolutions of their tangent cones. The Betti sequence of the resolutions of the tangent cones of these families is $(1,10,23,20,6)$.

\end{document}